\documentclass[leqno,12pt]{amsart}
\usepackage{amssymb} 
\theoremstyle{plain} 
\newtheorem{theorem}{\indent\sc Theorem}[section] 

\newtheorem{corollary}[theorem]{\indent\sc Corollary}
\newtheorem{proposition}[theorem]{\indent\sc Proposition}

\theoremstyle{definition} 

\newtheorem{remark}[theorem]{\indent\sc Remark}

\begin{document}

\title[Coupled Painlev\'e II systems]{On some Hamiltonian structures of coupled Painlev\'e II systems in dimension four \\}

\renewcommand{\thefootnote}{\fnsymbol{footnote}}
\footnote[0]{2000\textit{ Mathematics Subjet Classification}.
 14E05,14E15,20F55,34M55.}

\maketitle

\begin{abstract}
We find and study a two-parameter family of coupled Painlev\'e II systems in dimension four with affine Weyl group symmetry of several types. Moreover, we find a three-parameter family of polynomial Hamiltonian systems in two variables $t,s$. Setting $s=0$, we can obtain an autonomous version of the coupled Painlev\'e II systems. We also show its symmetry and holomorphy conditions.
\end{abstract}

\section{ Introduction}
In this paper, we study coupled Painlev\'e II systems in dimension four
\begin{equation}
\frac{dx}{dt}=\frac{\partial H}{\partial y}, \ \  \frac{dy}{dt}=-\frac{\partial H}{\partial x}, \ \  \frac{dz}{dt}=\frac{\partial H}{\partial w}, \ \  \frac{dw}{dt}=-\frac{\partial H}{\partial z}
\end{equation}
with the polynomial Hamiltonian
\begin{equation}
H=H_{II}(x,y,t;\alpha_1)+H_{II}(z,w,t;\alpha_2)+ayw.
\end{equation}
Here $x,y,z$ and $w$ denote unknown complex variables and $\alpha_1,\alpha_2$ are complex parameters and a is a coupling constant. The symbol $H_{II}(q,p,t;\alpha)$ denotes the Hamiltonian of the second-order Painlev\'e II systems given by
\begin{equation}
H_{II}(x,y,t;\alpha)=a_1x^2y+a_2y^2/2+a_3ty+a_1\alpha x,
\end{equation}
where $\alpha$ is a parameter and $a_i \ (i=1,2,3)$ are nonzero parameters which can be fixed arbitrarily.

As an important fact, for any $a \in {\Bbb C}$ the system (1) passes the Painlev\'e test. Moreover, the system (1) becomes again a polynomial Hamiltonian system in each coordinate system $(x_i,y_i,z_i,w_i) \ (i=1,2)$:
$$
(x_1,y_1,z_1,w_1)=(1/x,-(xy+\alpha_1)x,z,w), \ (x_2,y_2,z_2,w_2)=(x,y,1/z,-(zw+\alpha_2)z).
$$
Each coordinate system contains a three-parameter family of meromorphic solutions of the system (1).

Moreover, the system (1) is invariant under the following birational and symplectic transformations:
\begin{align}
\begin{split}
&S_1: (x,y,z,w,t;\alpha_1,\alpha_2) \rightarrow (x+\frac{\alpha_1}{y},y,z,w,t;-\alpha_1,\alpha_2),\\
&S_2: (x,y,z,w,t;\alpha_1,\alpha_2) \rightarrow (x,y,z+\frac{\alpha_2}{w},w,t;\alpha_1,-\alpha_2). 
\end{split}
\end{align}

In the case of $a=3/4$, A.N.W. Hone presented coupled Painlev\'e II systems with affine Weyl group symmetry of type $C_2^{(1)}$ (see \cite{Ho}). This system is explicitly written as
\begin{equation*}
  \left\{
  \begin{aligned}
   \frac{dx}{dt} &=\frac{\partial H}{\partial y}=2x^2+y/4+3w/4-t,\\
   \frac{dy}{dt} &=-\frac{\partial H}{\partial x}=-4xy-\alpha_1,\\
   \frac{dz}{dt} &=\frac{\partial H}{\partial w}=2z^2+w/4+3y/4-t,\\
   \frac{dw}{dt} &=-\frac{\partial H}{\partial z}=-4zw-\alpha_2
   \end{aligned}
  \right. 
\end{equation*}
wiht the polynomial Hamiltonian
\begin{equation}
H=2x^2y+y^2/8-ty+\alpha_1x+2z^2w+w^2/8-tw+\alpha_2z+3yw/4.
\end{equation}

In the case of $a=1$, in 2005, the author found a 2-parameter family of coupled Painlev\'e II systems with affine Weyl group symmetry of type $C_2^{(1)}$ (see \cite{Sasa1}). This system is explicitly written as
\begin{equation*}
  \left\{
  \begin{aligned}
   \frac{dx}{dt} &=\frac{\partial H}{\partial y}=-x^2+y+w-t/2,\\
   \frac{dy}{dt} &=-\frac{\partial H}{\partial x}=2xy+\alpha_1,\\
   \frac{dz}{dt} &=\frac{\partial H}{\partial w}=-z^2+y+w-t/2,\\
   \frac{dw}{dt} &=-\frac{\partial H}{\partial z}=2zw+\alpha_2
   \end{aligned}
  \right. 
\end{equation*}
wiht the polynomial Hamiltonian
\begin{equation}
H=-x^2y+y^2/2-ty/2-\alpha_1x-z^2w+w^2/2-tw/2-\alpha_2z+yw.
\end{equation}

In this paper, we study the case of $a=-3$. This system is explicitly given by
\begin{equation}\label{sys:1}
  \left\{
  \begin{aligned}
   \frac{dx}{dt} &=\frac{\partial H}{\partial y}=-2x^2+4y-3w-2t,\\
   \frac{dy}{dt} &=-\frac{\partial H}{\partial x}=4xy+2\alpha_2,\\
   \frac{dz}{dt} &=\frac{\partial H}{\partial w}=z^2+2w-3y+t,\\
   \frac{dw}{dt} &=-\frac{\partial H}{\partial z}=-2zw-\alpha_3
   \end{aligned}
  \right. 
\end{equation}
with the polynomial Hamiltonian
\begin{equation}\label{PH}
H=-2x^2y+2y^2-2ty-2\alpha_2x+z^2w+w^2+tw+\alpha_3z-3yw.
\end{equation}

This paper is organized as follows. In Section 1, we study symmetry of the system \eqref{sys:1}. In Section 2, we will study holomorphy condition of the system \eqref{sys:1}. In Section 3, we study some Hamiltonian structures of the system \eqref{sys:1}. In final section, we find a 3-parameter family of polynomial Hamiltonian systems in two variables $t,s$. Setting $s=0$, we can obtain an autonomous version of the system \eqref{sys:1}. We also show its symmetry and holomorphy conditions.

\section{Symmetry of the system \eqref{sys:1}}

In this section, we study symmetry of the system \eqref{sys:1}.

\begin{theorem}\label{th:1}
The system \eqref{sys:1} is invariant under the following birational and symplectic transformations{\rm : \rm}
\begin{align}
\begin{split}
s_1: &(x,y,z,w,t;\alpha_1,\alpha_2,\alpha_3) \rightarrow \left(x+\frac{\alpha_2}{y},y,z,w,t;\alpha_1+2\alpha_2,-\alpha_2,\alpha_3 \right),\\
s_2: &(x,y,z,w,t;\alpha_1,\alpha_2,\alpha_3) \rightarrow \left(x,y,z+\frac{\alpha_3}{w},w,t;\alpha_1+\alpha_3,\alpha_2,-\alpha_3 \right),\\
s_3: &(x,y,z,w,t;\alpha_1,\alpha_2,\alpha_3) \rightarrow\\
&(\frac{x(y-x^2-w-t)-(x+z)w+\alpha_1}{y-x^2-w-t},\\
&y-x^2-w+\frac{(x(y-x^2-w-t)-(x+z)w+\alpha_1)^2}{(y-x^2-w-t)^2}\\
&-\frac{(x+z)(-1-(x+z)w+2\alpha_1+2\alpha_2)}{y-x^2-w-t},\\
&-\frac{y-x^2-w-t}{x+z}-\frac{x(y-x^2-w-t)-(x+z)w+\alpha_1}{y-x^2-w-t},\\
&-\frac{(x+z)(-1-(x+z)w+2\alpha_1+2\alpha_2)}{y-x^2-w-t},t;-\alpha_1-\alpha_3,2\alpha_1+\alpha_2+\alpha_3,\alpha_3).
\end{split}
\end{align}
\end{theorem}
We note that the parameters $\alpha_1,\alpha_2,\alpha_3$ satisfy the following relation:
$$
2\alpha_1+2\alpha_2+\alpha_3=1.
$$
Theorem \ref{th:1} can be checked by a direct calculation.

\begin{corollary}
The transformations described in Theorem \ref{th:1} satisfy the following relations{\rm : \rm}
$$
s_1^2=s_2^2=s_3^2=(s_1s_2)^2=1.
$$
Moreover, the transformation $s_3s_1$ acts on parameters $(\alpha_1,\alpha_2,\alpha_3)$ as follows{\rm : \rm}
$$
s_3s_1:(\alpha_1,\alpha_2,\alpha_3) \rightarrow (\alpha_1+1,\alpha_2-1,\alpha_3).
$$
\end{corollary}

\begin{proposition}
The system \eqref{sys:1} has the following invariant cycles\rm{:\rm}
\begin{center}
\begin{tabular}{|c|c|c|} \hline
codimension & invariant cycle & parameter's relation \\ \hline
1 & $f_1:=y$ & $\alpha_2=0$  \\ \hline
1 & $f_2:=w$ & $\alpha_3=0$  \\ \hline
2 & $f_3^{(1)}:=y-x^2-w-t, \ f_3^{(2)}:=x+z$ & $\alpha_1=0$  \\ \hline
\end{tabular}
\end{center}
\end{proposition}

\section{Holomorphy of the system \eqref{sys:1}}

\begin{theorem}\label{th:2}
Let us consider a Hamiltonian system with polynomial Hamiltonian $H \in {\Bbb C}(t)[x,y,z,w]$. We assume that

$(A1)$ $deg(H)=5$ with respect to $x,y,z,w$.

$(A2)$ This system becomes again a polynomial Hamiltonian system in each coordinate system $(x_i,y_i,z_i,w_i) \ (i=1,2,3)${\rm : \rm}
$$
x_1=\frac{1}{x}, \ y_1=-x(xy+\alpha_2), \ z_1=z, \ w_1=w,
$$
$$
x_2=x,\ y_2=y,\ z_2=\frac{1}{z},\ w_2=-(zw+\alpha_3)z,
$$
$$
x_3=\frac{1}{x}, \ y_3=-((y-x^2-w-t)x-(x+z)w+\alpha_1)x, \ z_3=-\frac{w}{x}, \ w_3=x(x+z).
$$
Then such a system coincides with the system \eqref{sys:1} with the polynomial Hamiltonian \eqref{PH}.
\end{theorem}

By solving a linear problem, Theorem \ref{th:2} can be checked by a direct calculation.

\begin{remark}
\rm{Each coordinate system given in Theorem \ref{th:2} contains a 3-parameter family of meromorphic solutions of \eqref{sys:1}. \rm}
\end{remark}

\begin{proposition} 
On each affine open set $(x_j,y_j,z_j,w_j) \in U_j \times B$ in Theorem \ref{th:2}, each  Hamiltonian $H_j$ on $U_j \times B$ is expressed as a polynomial in $x_j,y_j,z_j,w_j,t$, and satisfies the following condition{\rm: \rm}
$$
dx \wedge dy +dz \wedge dw - dH \wedge dt=dx_j \wedge dy_j +dz_j \wedge dw_j - dH_j \wedge dt \ \ \ \ (j=1,2),
$$
$$
dx \wedge dy +dz \wedge dw - d(H+x) \wedge dt=dx_3 \wedge dy_3 +dz_3 \wedge dw_3 - dH_3 \wedge dt.
$$
\end{proposition}

\section{Hamiltonian structures}
There is the following birational and symplectic transformations other than the transformations given in Theorems \ref{th:1} and \ref{th:2}.

\begin{theorem}\label{th:3}
The system \eqref{sys:1} is equivalent to the following Hamiltonian systems{\rm : \rm}

by using the rational and symplectic transformation $\varphi_1$
\begin{align}
\begin{split}
\varphi_1:(x,y,z,w) \longrightarrow &(x-\frac{(x+z)w-\alpha_1}{-t-w-x^2+y},-t-w-x^2+y,\\
&-w(-t-w-x^2+y),\frac{x+z}{-t-w-x^2+y}),
\end{split}
\end{align}
the system \eqref{sys:1} is transformed to
\begin{equation}
  \left\{
  \begin{aligned}
   \frac{dx}{dt} &=2x^2+4y+2t-zw^2+\alpha_3w,\\
   \frac{dy}{dt} &=-4xy+2zw+2\alpha_1,\\
   \frac{dz}{dt} &=-2xz-2yzw+\alpha_3y,\\
   \frac{dw}{dt} &=yw^2+2xw+1
   \end{aligned}
  \right. 
\end{equation}
with the polynomial Hamiltonian
\begin{equation}
H=2x^2y+2y^2+2ty-yzw^2+\alpha_3yw-2xzw-2\alpha_1x-z,
\end{equation}
and by using the symplectic transformation $\varphi_2$
\begin{equation}
\varphi_2:(x,y,z,w) \longrightarrow (x,-t-w-x^2+y,-w,x+z),
\end{equation}
the system \eqref{sys:1} is transformed to
\begin{equation}\label{sys:4}
  \left\{
  \begin{aligned}
   \frac{dx}{dt} &=2x^2+4y-z+2t,\\
   \frac{dy}{dt} &=-4xy-2zw-2\alpha_1,\\
   \frac{dz}{dt} &=-2zw+2xz+\alpha_3,\\
   \frac{dw}{dt} &=w^2-2xw+y
   \end{aligned}
  \right. 
\end{equation}
with the polynomial Hamiltonian
\begin{equation}
H=2x^2y+2y^2-yz+2ty+2\alpha_1x+2xzw+\alpha_3w-zw^2.
\end{equation}
\end{theorem}
If $\alpha_1=0$, the system \eqref{sys:4} has a particular solution $y=w=0$. Here the systen in the variables $(x,z)$ satisfies
\begin{equation}
  \left\{
  \begin{aligned}
   \frac{dx}{dt} &=2x^2-z+2t,\\
   \frac{dz}{dt} &=2xz+\alpha_3.
   \end{aligned}
  \right. 
\end{equation}

Theorem \ref{th:3} can be checked by a direct calculation.

\section{Autonomous version of the system \eqref{sys:1} and partial differential system in two variables}

In this section, we find a 3-parameter family of partial differential systems in two variables $t,s$ given by
\begin{equation}\label{eq:10}
  \left\{
  \begin{aligned}
   dq_1 =&\frac{\partial K_1}{\partial p_1}dt+\frac{\partial K_2}{\partial p_1}ds,\\
   dp_1 =&-\frac{\partial K_1}{\partial q_1}dt-\frac{\partial K_2}{\partial q_1}ds,\\
   dq_2 =&\frac{\partial K_1}{\partial p_2}dt+\frac{\partial K_2}{\partial p_2}ds,\\
   dp_2 =&-\frac{\partial K_1}{\partial q_2}dt-\frac{\partial K_2}{\partial q_2}ds
   \end{aligned}
  \right. 
\end{equation}
with the polynomial Hamiltonians
\begin{align}
\begin{split}
K_1=&q_1^2p_1-p_1^2+\alpha_2 q_1-\frac{q_2^2 p_2}{2}-\frac{p_2^2}{2}-\frac{\alpha_3}{2}q_2+\frac{3}{2}p_1 p_2,\\
K_2=&-\alpha_3^2 p_1+q_2^4p_2^2+2q_2^2p_2^3+p_2^4+2\alpha_3 q_2^3 p_2+4(\alpha_1+\alpha_3)q_2 p_2^2+\alpha_3^2 q_2^2+\alpha_3(2\alpha_1+\alpha_3)p_2\\
&-p_2(8q_1p_1q_2p_2+6p_1q_2^2p_2+4q_1^2p_1p_2+2p_1p_2^2-p_1^2p_2+4\alpha_2 q_1p_2+4\alpha_3 q_1p_1+6\alpha_3 p_1q_2).
\end{split}
\end{align}
Setting $s=0$, we can obtain an autonomous version of the system \eqref{sys:1}.

\begin{proposition}
The system \eqref{eq:10} satisfies the compatibility conditions$:$
\begin{equation}
\frac{\partial }{\partial s} \frac{\partial q_1}{\partial t}=\frac{\partial }{\partial t} \frac{\partial q_1}{\partial s}, \quad \frac{\partial }{\partial s} \frac{\partial p_1}{\partial t}=\frac{\partial }{\partial t} \frac{\partial p_1}{\partial s}, \quad \frac{\partial }{\partial s} \frac{\partial q_2}{\partial t}=\frac{\partial }{\partial t} \frac{\partial q_2}{\partial s}, \quad \frac{\partial }{\partial s} \frac{\partial p_2}{\partial t}=\frac{\partial }{\partial t} \frac{\partial p_2}{\partial s}.
\end{equation}
\end{proposition}

\begin{proposition}
The system \eqref{eq:10} has $K_1$ and $K_2$ as its first integrals.
\end{proposition}

\begin{proposition}
Two Hamiltonians $K_1$ and $K_2$ satisfy
\begin{equation}
\{K_1,K_2\}=0,
\end{equation}
where
\begin{equation}
\{K_1,K_2\}=\frac{\partial K_1}{\partial p_1}\frac{\partial K_2}{\partial q_1}-\frac{\partial K_1}{\partial q_1}\frac{\partial K_2}{\partial p_1}+\frac{\partial K_1}{\partial p_2}\frac{\partial K_2}{\partial q_2}-\frac{\partial K_1}{\partial q_2}\frac{\partial K_2}{\partial p_2}.
\end{equation}
\end{proposition}
Here, $\{,\}$ denotes the poisson bracket such that $\{p_i,q_j\}={\delta}_{ij}$ (${\delta}_{ij}$:kronecker's delta).

\begin{proposition}
The system \eqref{eq:10} is invariant under the following birational and symplectic transformations{\rm : \rm} with {\it the notation} $(*):=(q_1,p_1,q_2,p_2,t,s;\alpha_1,\alpha_2,\alpha_3)$\rm{: \rm}
\begin{align}
\begin{split}
s_1: &(*) \rightarrow \left(q_1+\frac{\alpha_2}{p_1},p_1,q_2,p_2,t,s;\alpha_1+2\alpha_2,-\alpha_2,\alpha_3 \right),\\
s_2: &(*) \rightarrow \left(q_1,p_1,q_2+\frac{\alpha_3}{p_2},p_2,t,s;\alpha_1+\alpha_3,\alpha_2,-\alpha_3 \right),\\
s_3: &(*) \rightarrow (\frac{q_1(p_1-q_1^2-p_2)-(q_1+q_2)p_2+\alpha_1}{p_1-q_1^2-p_2},\\
&p_1-q_1^2-p_2+\frac{(q_1(p_1-q_1^2-p_2)-(q_1+q_2)p_2+\alpha_1)^2}{(p_1-q_1^2-p_2)^2}\\
&-\frac{(q_1+q_2)(-(q_1+q_2)p_2+2\alpha_1+2\alpha_2)}{p_1-q_1^2-p_2},\\
&-\frac{p_1-q_1^2-p_2}{q_1+q_2}-\frac{q_1(p_1-q_1^2-p_2)-(q_1+q_2)p_2+\alpha_1}{p_1-q_1^2-p_2},\\
&-\frac{(q_1+q_2)(-(q_1+q_2)p_2+2\alpha_1+2\alpha_2)}{p_1-q_1^2-p_2},t,s;-\alpha_1-\alpha_3,2\alpha_1+\alpha_2+\alpha_3,\alpha_3).
\end{split}
\end{align}
\end{proposition}
Here, the parameters $\alpha_i$ satisfy the relation $2\alpha_1+2\alpha_2+\alpha_3=0$.

\begin{corollary}
The transformation $s_3s_1$ acts on parameters $(\alpha_1,\alpha_2,\alpha_3)$ as follows{\rm : \rm}
$$
s_3s_1:(\alpha_1,\alpha_2,\alpha_3) \rightarrow (\alpha_1,\alpha_2,\alpha_3).
$$
\end{corollary}
We see that this transformation is an automorphism of the system \eqref{eq:10}. Applying these B{\"a}cklund transformations $(s_3s_1)^m \ (m=1,2,\ldots)$, we can obtain a series of its particular solutions.

\begin{proposition}
Let us consider a polynomial Hamiltonian system with Hamiltonian $K \in {\Bbb C}[q_1,p_1,q_2,p_2]$. We assume that

$(B1)$ $deg(K)=6$ with respect to $q_1,p_1,q_2,p_2$.

$(B2)$ This system becomes again a polynomial Hamiltonian system in each coordinate $R_i \ (i=1,2,3)${\rm : \rm}
\begin{align*}
\begin{split}
R_1:(x_1,y_1,z_1,w_1)=&\left(\frac{1}{q_1},-(q_1p_1+\alpha_2)q_1,q_2,p_2 \right),\\
R_2:(x_2,y_2,z_2,w_2)=&\left(q_1,p_1,\frac{1}{q_2},-(q_2p_2+\alpha_3)q_2 \right),\\
R_3:(x_3,y_3,z_3,w_3)=&\left(\frac{1}{q_1},-((p_1-q_1^2-p_2)q_1-(q_1+q_2)p_2+\alpha_1)q_1,-\frac{p_2}{q_1},q_1(q_1+q_2) \right),
\end{split}
\end{align*}
where the parameters $\alpha_i$ satisfy the relation $2\alpha_1+2\alpha_2+\alpha_3=0$. Then such a system coincides with the Hamiltonian system \eqref{eq:10} with the polynomial Hamiltonians $K_1,K_2$.
\end{proposition}
We remark that we can obtain three polynomial Hamiltonians $K_1,K_2$ and $K_3$ satisfying this proposition, where the Hamiltonian $K_3$ satisfy the relation
$$
K_3=\frac{1}{4}(4K_1^2-13K_2).
$$
We note that the conditions $(B2)$ should be read that
\begin{align*}
\begin{split}
&R_i(K) \quad (i=1,2,3)
\end{split}
\end{align*}
are polynomials with respect to $x_i,y_i,z_i,w_i$.


\begin{thebibliography}{99}
\bibitem[1]{1} P. Painlev\'e, {\em M\'emoire sur les \'equations diff\'erentielles dont l'int\'egrale g\'en\'erale est uniforme}, Bull. Soci\'et\'e Math\'ematique de France. {\bf 28} (1900),  201--261.

\bibitem[2]{2} P. Painlev\'e, {\em Sur les \'equations diff\'erentielles du second ordre et d'ordre sup\'erieur dont l'int\'egrale est uniforme}, Acta Math. {\bf 25} (1902), 1--85. 

\bibitem[3]{3} B. Gambier, {\em Sur les \'equations diff\'erentielles du second ordre et du premier degr\'e dont l'int\'egrale g\'en\'erale est \`a points critiques fixes}, Acta Math. {\bf 33} (1910), 1--55.


\bibitem[4]{Cosgrove1} C. M. Cosgrove and G. Scoufis,
{\em Painlev\'e classification of a class of differential equations of the second order and second degree}, Studies in Applied Mathematics. {\bf 88} (1993), 25-87.

\bibitem[5]{Cosgrove2} C. M. Cosgrove,
{\em All binomial-type Painlev\'e equations of the second order and degree three or higher}, Studies in Applied Mathematics. {\bf 90} (1993), 119-187.

\bibitem[6]{6} F. Bureau, 
{\em Integration of some nonlinear systems of ordinary differential equations}, 
Annali di Matematica. {\bf 94} (1972), 345--359. 

\bibitem[7]{7} J. Chazy, 
{\em Sur les \'equations diff\'erentielles dont l'int\'egrale g\'en\'erale est uniforme et admet des singularit\'es essentielles mobiles}, 
Comptes Rendus de l'Acad\'emie des Sciences, Paris. {\bf 149} (1909), 563--565. 

\bibitem[8]{8} J. Chazy, 
{\em Sur les \'equations diff\'erentielles dont l'int\'egrale g\'en\'erale poss\'ede une coupure essentielle mobile }, 
Comptes Rendus de l'Acad\'emie des Sciences, Paris. {\bf 150} (1910), 456--458. 


\bibitem[9]{9} J. Chazy, 
{\em Sur les \'equations diff\'erentielles du trousi\'eme ordre et d'ordre sup\'erieur dont l'int\'egrale a ses points critiques fixes}, 
Acta Math. {\bf 34} (1911), 317--385. 



\bibitem[10]{Sasa0} Y. Sasano, {\em Coupled Painlev\'e VI systems in dimension four with affine Weyl group symmetry of types $B_6^{(1)},D_6^{(1)}$ and $D_7^{(2)}$}, preprint.

\bibitem[11]{Sasa2} Y. Sasano, {\em Four-dimensional Painlev\'e systems of types $D_5^{(1)}$ and $B_4^{(1)}$}, preprint.

\bibitem[12]{Sasa5} Y. Sasano, {\em Higher order Painlev\'e equations of type ${D_l}^{(1)}$}, RIMS Kokyuroku {\bf 1473} (2006), 143--163.

\bibitem[13]{Sasa4} Y. Sasano, {\em Symmetries in the system of type $D_4^{(1)}$}, preprint.

\bibitem[14]{Sasa3} Y. Sasano, {\em Coupled Painlev\'e III systems with affine Weyl group symmetry of types $B_4^{(1)}$, $D_4^{(1)}$ and $D_5^{(2)}$}, preprint.

\bibitem[15]{Sasa6} Y. Sasano, {\em Coupled Painlev\'e III systems with affine Weyl group symmetry of types $B_5^{(1)},D_5^{(1)}$ and $D_6^{(2)}$}, preprint.

\bibitem[16]{Sasa7} Y. Sasano, {\em Coupled Painlev\'e VI systems in dimension four with affine Weyl group symmetry of type $D_6^{(1)}$, II}, RIMS Kokyuroku Bessatsu. {\bf B5} (2008), 137--152.

\bibitem[17]{Sasa8} Y. Sasano, {\em Coupled Painlev\'e VI systems in dimension four with affine Weyl group symmetry of type $E_6^{(2)}$}, preprint.



\bibitem[18]{Sasa10} Y. Sasano, {\em Symmetries in the system of type $A_5^{(2)}$}, preprint.

\bibitem[19]{Sasa11} Y. Sasano, {\em Coupled Painlev\'e systems with affine Weyl group symmetry of types $A_7^{(2)},A_5^{(2)}$ and $D_4^{(3)}$}, preprint.

\bibitem[20]{Sasa12} Y. Sasano, {\em Symmetry and holomorphy of the second member of the second Painlev\'e hierarchy}, preprint.


\bibitem[21]{Sasa1} Y. Sasano, {\em Coupled Painlev\'e II system in dimension four and the systems of type $A_4^{(1)}$}, Tohoku Math J. {\bf 58} (2006), 529--548.

\bibitem[22]{Ho} A.N.W. Hone, {\em Coupled Painlev\'e systems and quartic potentials}, J.Phys.A. Gen. {\bf 34} (2001), 2235--2245.


\end{thebibliography}
\end{document}